\pdfoutput=1  

\documentclass{article}
\usepackage{graphicx}
\usepackage{wrapfig}
\usepackage{subcaption}
\usepackage[T1]{fontenc}
\usepackage{helvet}
\usepackage{mathbbol}
\usepackage{mathtools}

\usepackage{booktabs}
\usepackage{multirow}
\usepackage{hyperref}
\usepackage{float}
\usepackage{boldline} 
\usepackage{hhline}
\usepackage{longtable}
\usepackage{iclr2026_conference, times}
\iclrfinalcopy

\usepackage[maxbibnames=1, giveninits=true, sorting=nyt, style=authoryear, maxcitenames=1, uniquelist=false,
  mincitenames=1]{biblatex}
\AtEveryBibitem{
  \clearfield{doi}
  \clearfield{isbn}
  \clearfield{publisher}
  \clearfield{address}
}

\usepackage{enumitem}

\usepackage[most]{tcolorbox}

\DeclareMathOperator{\Li}{Li}

\title{The Ramanujan Challenge For AI}

\author{
\\[-20pt]
\textbf{Ramanujan Machine Group\textsuperscript{*}}
}

\begin{document}

\maketitle
\vspace{-20pt}

\begin{abstract}
\vspace{-6pt}

To help evaluate the mathematical skills of current AI systems, we present a set of formulas for fundamental mathematical constants. These problems are attractive for AI evaluation because they are concrete and can be checked numerically to arbitrary precision, yet proving them may require non-obvious mathematics. Mathematical constants such as $\pi$, $e$, Catalan’s constant, and special values of the Riemann zeta function have fascinated mathematicians for centuries. The search for formulas evaluating mathematical constants has produced some of the most beautiful mathematics in the field, especially in cases that yield irrationality proofs or fast convergence rates. Ramanujan’s legacy is emblematic of this tradition. The list we provide contains two types of problems: formulas whose proofs are known to the authors but will remain encrypted for a short initial period; and formulas that are not yet proven. We are curious to see the achievements of AI in both cases.

\end{abstract}

\vspace{-8pt}

\begin{center}
\renewcommand{\arraystretch}{1.2}
\begin{tabular}{c p{0.59\textwidth} p{0.21\textwidth}}
\toprule[1pt]
\multicolumn{1}{c}{\textbf{Problem}} &
\multicolumn{1}{c}{\textbf{Name}} &
\multicolumn{1}{c}{\textbf{Contributor}} \\
\midrule[0.8pt] 

2.1 & Polynomial continued fraction for $\pi$ 
& Michael Shalyt\textsuperscript{1} \\

\midrule[0.2pt] 

2.2 & Euler's constant $\gamma$ as an Ap\'ery limit 
& \multirow{2}{0.21\textwidth}{\raggedright Rotem Kalisch\textsuperscript{1}} \\

2.3 & The sum $\pi+e$ as an Ap\'ery limit 
& \\

\midrule[0.2pt]

2.4 & A series of harmonic numbers converging to a polylogarithm combined with zeta values 
& \multirow{2}{0.21\textwidth}{\raggedright Carsten Schneider\textsuperscript{2}}\\

\midrule[0.2pt]

2.5 & Efficient rational approximation of Catalan's constant $G$ 
& \multirow{2}{0.21\textwidth}{\raggedright Hila Barkan\textsuperscript{1}} \\

2.6 & A series for $\zeta(2)+\zeta(3)$ 
& \\

\midrule[0.2pt]

2.7 & Efficient four-term recurrence for $\zeta(2)+\zeta(3)$ 
& \multirow{2}{0.21\textwidth}{\raggedright Elyasheev Leibtag\textsuperscript{3}} \\

2.8 & Very fast rational approximation of $\sqrt{10005}/\pi$ 
& \\

\midrule[0.2pt]

3.1 & An integral over knot polynomial roots expressing $\pi^2$ 
& John Campbell\textsuperscript{4} \\

\midrule[0.2pt]

3.2 & Optimality of Ap\'ery's irrationality-measure bound for $\zeta(3)$ 
& Shachar Weinbaum\textsuperscript{1} \\

\bottomrule[1pt]
\end{tabular}
\end{center}

\vspace{-8pt}

\begin{center}
\renewcommand{\arraystretch}{1.2}
\begin{tabular}{p{0.25\textwidth} p{0.25\textwidth} p{0.32\textwidth}}
\toprule
\multicolumn{1}{c}{\textbf{Role}} &
\multicolumn{1}{c}{\textbf{Name}} &
\multicolumn{1}{c}{\textbf{Email}}\\
\midrule
Proof Collection & Tali Monderer\textsuperscript{1} & talimon@campus.technion.ac.il \\
Validation & Ashvni Narayanan\textsuperscript{1} & ashvni.n@campus.technion.ac.il \\
Principal Investigator & Ido Kaminer\textsuperscript{1} & kaminer@technion.ac.il \\
\bottomrule
\end{tabular}
\end{center}

\vspace{-6pt}

{\footnotesize
\textsuperscript{1}Ramanujan Machine Group, Technion, Haifa, Israel. \\
\textsuperscript{2}Research Institute for Symbolic Computation, Johannes Kepler Universit\"at Linz, Linz, Austria. \\
\textsuperscript{3}Department of Mathematics, Computer Science and Statistics, Ghent University, Ghent, Belgium.
\\
\textsuperscript{4}Department of Mathematics and Statistics, Toronto Metropolitan University, Toronto, Canada.
\\
\textsuperscript{*}ramanujan.machine@gmail.com
}

\section{Introduction}
\label{section-introduction}

Recent progress in AI for mathematics has made quantitative evaluation increasingly urgent. We increasingly need university-level and research-level questions that test whether AI systems can contribute to genuine mathematical work.
Several recent benchmarks address this need from different directions. 
RealMath and LemmaBench study mathematical questions drawn from research papers and mathematical forums [\cite{ZhangEtAl2025, peyronnet2026lemmabenchliveresearchlevelbenchmark}].
A central difficulty is contamination: if a problem or its solution appears in the AI training data, success may reflect retrieval rather than reasoning. One response is to perturb existing problems, as in GSM-Symbolic [\cite{GSM-Symbolic-Numeric-Variation}] and ASyMOB [\cite{ASyMOB}], reducing dependence on uncontaminated original questions. 
A stronger response is to use authentic research problems that have not yet appeared publicly. 
FrontierMath [\cite{GlazerEtAl2024}], Riemann-Bench [\cite{garre2026riemannbenchbenchmarkmoonshotmathematics}], and part of Humanity's Last Exam [\cite{Humanity-Last_Exam}] emphasize original difficult questions with structured verification, while First Proof [\cite{AbouzaidEtAl2026}] focuses on unpublished research problems whose answers were known to experts.

The present manuscript follows the spirit of First Proof, but concentrates on a more focused domain with a long mathematical tradition: explicit formulas involving fundamental mathematical constants, presented as recurrences, continued fractions, series, and integrals.
The most interesting formulas often point to hidden structures (algebraic, group-theoretic, etc.).
Historically, such formulas have been signatures of mathematical inspiration, from Ramanujan’s ``dreams'' to Apéry’s ``garden''. 

This domain is well suited to quantitative evaluation of AI abilities. On the one hand, candidate formulas for constants can often be tested to thousands of digits, so numerical validation is immediate and objective. 
On the other hand, converting such evidence into proof may require a wide variety of mathematical tools from various domains. It is often hard to know in advance which approach will be best suited for each formula.

The historical examples are iconic. Ramanujan's formulas for $\pi$ remain a model of unexpected structure in special values [\cite{Ramanujan1914}]. Ap\'ery's proof of the irrationality of $\zeta(3)$, together with the reinterpretation by Beukers via multiple integrals, showed how recurrence relations and arithmetic estimates can turn an experimental pattern into a theorem [\cite{Apery1979,Beukers1979}]. Subsequent work and continued ongoing efforts extend this perspective to odd zeta values, Catalan's constant, Euler's constant, and a wide range of Ap\'ery limits and related special values [\cite{BallRivoal2001,Zudilin2001,RivoalZudilin2003,Lagarias2013,ChamberlandStraub2021,Raayoni2021,elimelech2023algorithm,weinbaum2025conservative}].

Many of the questions collected here lie in or near the modern Ap\'ery-style tradition. Some ask for polynomial continued fractions or holonomic recurrences converging to constants such as $\pi$, $\log 2$, $\gamma$, $\pi+\mathrm{e}$, or combinations of zeta values. Others concern efficient recurrences with a conjectural exponential rate of convergence, or integrals whose closed forms appear to encode deeper arithmetic structure. 
At a broader conceptual level, these identities also belong to the world of periods and special values, where numerical discovery frequently precedes structural explanation [\cite{KontsevichZagier2001}].

\subsection*{A few remarks}
The present list is split into two categories. The first category consists of formulas for which a proof (in the classic sense or as the result of an algorithmic procedure, see discussion below) is  known to the authors but has not yet been made public. These questions are intended to test whether an AI system, or a human participant, can independently structure a valid proof. The second category consists of formulas and hypotheses that have been numerically validated to high precision, but for which the authors do not currently have a proof. These are posed as open problems.

\textbf{Why are the formulas so complex?} Mostly because AI got so good, 
and can immediately prove most simple forms.
Nevertheless, some of the most interesting formulas presented below are complex because no simpler forms are known, as in the example of the recurrence that produces $\pi+e$. We will see whether/when AI can provide new insight into these classic questions.

\section{The Questions}
\label{section-questions}

This section contains proven problems, presented here for the first time. The proofs or computational procedures for the proofs are known to the authors but are not yet public.
Their level of difficulty varies. Some may be tackled using existing literature or folklore, while others appear to require additional novelty. In each case, the statement is written in the form of an explicit formula so that it can be checked numerically.

\subsection{Polynomial continued fraction for $\pi$}
Let:
\[
a_n = -220n^3 - 484n^2 - 301n - 42,
\qquad
b_n = 4n^2(2n+1)^2(5n-4)(5n+6).
\]

Prove:
\[
a_0 + \cfrac{b_1}{a_1 + \cfrac{b_2}{a_2 + \cfrac{b_3}{a_3 + \ddots}}}  \;=\; \frac{6}{3-\pi}.
\]

\subsection{Euler's constant $\gamma$ as an Ap\'ery limit}
For $n\geq 0$, define:
\begin{align*}
0 = &\left(- 8 n^{3} - 51 n^{2} - 105 n - 68\right) u_n \\
& + \left(24 n^{5} + 337 n^{4} + 1833 n^{3} + 4818 n^{2} + 6092 n + 2928\right) u_{n-1} \\
&- \left(n + 2\right) \left(n + 3\right) \left(24 n^{5} + 273 n^{4} + 1150 n^{3} + 2154 n^{2} + 1635 n + 268\right) u_{n-2} \\
& +\left(n + 1\right) \left(n + 2\right)^{4} \left(n + 3\right) \left(8 n^{3} + 75 n^{2} + 231 n + 232\right) u_{n-3}\\ 
\end{align*}

Let $p_n, q_n$ be two solutions of the recurrence, defined by the initial values:
\[p_{-3} = 0, \quad p_{-2} = 7,  \quad p_{-1} = 179\]
\[q_{-3} = 1, \quad q_{-2} = 12, \quad q_{-1} = 306\]

Prove:
\[
\lim_{n\to\infty}{\frac{p_n}{q_n}}=\gamma,
\]
where $\gamma$ is Euler's constant.

\subsection{The sum $\pi+e$ as an Ap\'ery limit}

For $n\geq 1,$ define  
\begin{align*}
0 = & \left(- n^{3} + 2 n^{2} + 7 n + 3\right) u_n \\
&+ \left(n + 2\right) \left(2 n^{4} + n^{3} - 26 n^{2} - 48 n - 19\right) u_{n-1} \\
&+ \left(n + 2\right) \left(n^{6} + 9 n^{5} + 8 n^{4} - 87 n^{3} - 249 n^{2} - 234 n - 68\right) u_{n-2} \\ 
&+ \left(n + 1\right)^{2} \left(n + 2\right) \left(2 n^{5} + 3 n^{4} - 13 n^{3} - 21 n^{2} + 4\right) u_{n-3} \\
& - n^{3} \left(n + 1\right)^{2} \left(n + 2\right) \left(n^{3} + n^{2} - 8 n - 11\right) u_{n-4}.
\end{align*}

Let $p_n, q_n$ be two solutions of the recurrence, defined by the initial values:
\[p_{-3} = 1, \quad p_{-2} = 1, \quad p_{-1} = 20, \quad p_0 = 296\]
\[q_{-3} = 1, \quad q_{-2} = 0, \quad q_{-1} = 4, \quad q_0 = 48\]

Prove: 
\[
\lim_{n\to\infty}\frac{p_n}{q_n} = \pi + e.
\]

\subsection{A series of harmonic numbers converging to a polylogarithm combined with zeta values}
Prove:
\begingroup
\small
\[
\sum_{m=0}^{\infty}\sum_{k=0}^{m} \frac{\binom{m}{k}^2 H_k^2}{(m+1)^2\binom{2m}{m}} = 
\]
\[
20 \Li_4\!\left(\frac{1}{2}\right) + \frac{5}{6}\log^4(2) + 10 \zeta(2) - \frac{65}{9}\zeta(2)^2 - 
 \log^2(2) (12 + 5 \zeta(2)) 
 + \frac{1}{2}\zeta(3) 
 + \log(2)\left(\frac{35}{2}\zeta(3) - 16\right),
\]
\endgroup

where $H_0=0$ and for $k\geq 1$, $H_k$ is the $k$-th harmonic number.

\subsection{Efficient rational approximation of Catalan's constant $G$}
For $n\geq 0$, let

\[
M(n)=\bigl(m_{ij}(n)\bigr)_{1\leq i,j\leq 3},
\]
where
\[
\begin{aligned}
m_{11}(n)&=(-2n-5)(n+3)^2(136n^4+1424n^3+5548n^2+9551n+6141),\\
m_{12}(n)&=384n^6+6384n^5+44168n^4+162698n^3+336377n^2+369933n+169011,\\
m_{13}(n)&=-480n^4-4980n^3-19210n^2-32690n-20730,\\
m_{21}(n)&=(n+2)^2(n+3)^2(4n+10)(48n^3+386n^2+1017n+879),\\
m_{22}(n)&=(n+2)^2(-272n^5-3848n^4-21732n^3-61184n^2-85761n-47808),\\
m_{23}(n)&=(n+2)^2(320n^3+2540n^2+6610n+5640),\\
m_{31}(n)&=(-4n-10)(n+2)^2(n+3)^2(32n^4+302n^3+1037n^2+1530n+813),\\
m_{32}(n)&=(n+2)^2(192n^6+2984n^5+19116n^4+64452n^3+120256n^2+117279n+46476),\\
m_{33}(n)&=(n+2)^2(-16n^5-408n^4-2912n^3-8884n^2-12254n-6240).
\end{aligned}
\]

For $N\geq 0$, define
\[
\mathcal M_N
=
M(0)M(1)\cdots M(N-1).
\]

Choose initial conditions
\[
A=
\begin{pmatrix}
30921 & -32972 & 8240\\
33750 & -36000 & 9000
\end{pmatrix}.
\]
Write
\[
A\mathcal M_N
=
\begin{pmatrix}
P_{N,1} & P_{N,2} & P_{N,3}\\
Q_{N,1} & Q_{N,2} & Q_{N,3}
\end{pmatrix}.
\]
Then, for each $j=1,2,3$, prove:
\[
\lim_{N\to\infty}
\frac{P_{N,j}}{Q_{N,j}}
=
G,
\]
where
\[
G=\sum_{k=0}^{\infty}\frac{(-1)^k}{(2k+1)^2}
\]
is Catalan's constant.

\subsection{A series for $\zeta(2)+\zeta(3)$}
Let $(u_n)_{n\geq 1}$ be the sequence defined by
\[
u_1=-\frac{93}{4480},
\qquad
u_2=-\frac{117}{14000},
\]
and, for $n\geq 3$, by the recurrence

\begin{align*}
0=&-2(n+3)^3(2n+5)(3n+5)u_n \\
&\quad +(n+2)^2(15n^3+85n^2+155n+93)u_{n-1}\\
&\quad -(n+1)^3(n+2)(3n+8)u_{n-2}.
\end{align*}

Prove:
\[
\frac{2077}{720}+\sum_{j=1}^{\infty}u_j
=
\zeta(2)+\zeta(3).
\]

\subsection{Efficient four-term recurrence for $\zeta(2)+\zeta(3)$}

For \(n\ge 0\), define

$A_n = 1024(2n+5)^4(2n+7)^3(2n+9)^3\bigl(946n^2+6407n+10860\bigr)$,\\
$B_n = 128(2n+7)^3(2n+9)^3\bigl(104060n^6+1745370n^5+12145238n^4+44886481n^3+92943995n^2+102256019n+46709052\bigr)$,\\
$C_n = 16(n+3)^4(2n+9)^3\bigl(3784n^5+57792n^4+351019n^3+1059230n^2+1587211n+944620\bigr)$,\\
$D_n = (n+3)^4(n+4)^6\bigl(946n^2+4515n+5399\bigr)$.

Let \((p_n)_{n\ge 0}\) and \((q_n)_{n\ge 0}\) be the two solutions of the recurrence
\[
u_{n+1} = \frac{B_n}{A_n}\,u_n-\frac{C_{n-1}}{A_{n-1}}\,u_{n-1}+\frac{D_{n-2}}{A_{n-2}}\,u_{n-2}, \qquad n\ge 2,
\]
with initial conditions
\[
p_0=-612218384750,\qquad
p_1=-\frac{9525021973931919}{18100},\qquad
p_2=-\frac{29561828382772029}{65380},
\]
and
\[
q_0=-215040420000,\qquad
q_1=-\frac{167282265043404}{905},\qquad
q_2=-\frac{964185327658080}{6071}.
\]

Prove:
\[
\lim_{n\to\infty}\frac{p_n}{q_n}=\zeta(2)+\zeta(3).
\]

\subsection{Very fast rational approximation of $\sqrt{10005}/\pi$}

Let
\[
R=151931373056001,\qquad u=2n+3,\qquad
w=u(3u-2)(3u+2).
\]
For \(n\geq 0\), let
\[
M(n)=
\begin{pmatrix}
\dfrac{a_1}{w} & \dfrac{a_2}{w} & \dfrac{a_3}{w} & \dfrac{144R(u-1)^2}{w}\\[0.8em]
-u^3 & -3u^2 & -3u & -1\\[0.6em]
\dfrac{b_1}{144R} & -\dfrac{b_2}{72R} & -\dfrac{b_3}{36R} & -\dfrac{2u+2R-7}{2R}\\[0.8em]
\dfrac{c_1}{288R^2} & \dfrac{c_2}{144R^2} & \dfrac{c_3}{72R^2} & \dfrac{c_4}{4R^2}
\end{pmatrix},
\]
where
\[
\begin{aligned}
a_1&=(144R-99)u^5-(288R-333)u^4+(144R-229)u^3-114u^2+40u+64,\\
a_2&=(432R-243)u^4-(864R-909)u^3+(432R-868)u^2-80u+272,\\
a_3&=(432R-153)u^3-(864R-648)u^2+(432R-860)u+360,
\end{aligned}
\]
\[
\begin{aligned}
b_1&=9u^4-(144R-63)u^3+158u^2+168u+64,\\
b_2&=36u^3+(216R-189)u^2-316u-168,\\
b_3&=54u^2+(108R-189)u-158,
\end{aligned}
\]
and
\[
\begin{aligned}
c_1={}&18u^5+(54R+45)u^4-(288R^2-378R+251)u^3\\
&\quad +(948R-1086)u^2+(1008R-1384)u+(384R-576),\\
c_2={}&(153R-72)u^4-(657R-702)u^3-(432R^2-1292R+1069)u^2\\
&\quad +(2064R-2508)u+(1072R-1512),\\
c_3={}&(180R-108)u^3-(891R-864)u^2-(216R^2-1450R+1385)u\\
&\quad +(1116R-1422),\\
c_4={}&(6R-4)u^2-(33R-32)u-(4R^2-58R-236337691420383).
\end{aligned}
\]

For \(N\geq 0\), define
\[
\mathcal M_N=M(0)M(1)\cdots M(N-1),
\]
with the convention that \(\mathcal M_0=I\).

Choose initial conditions
\[
A=
\begin{pmatrix}
A_1&A_2&A_3&A_4\\
B_1&B_2&B_3&B_4
\end{pmatrix},
\]
where
\[
\begin{aligned}
(A_1,A_2,A_3,A_4)
={}&(
37169305760442252761441,\,
111507917281327441564208,\\
&\quad
111507917281327599720129,\,
37169305760442410917362),
\end{aligned}
\]
and
\[
\begin{aligned}
(B_1,B_2,B_3,B_4)
={}&(
1167416361542639692320,\,
3502249084627896132160,\\
&\quad
3502249084627879697280,\,
1167416361542622723840).
\end{aligned}
\]

Write
\[
A\mathcal M_N
=
\begin{pmatrix}
P_{N,1}&P_{N,2}&P_{N,3}&P_{N,4}\\
Q_{N,1}&Q_{N,2}&Q_{N,3}&Q_{N,4}
\end{pmatrix}.
\]
Then, for each \(j=1,2,3,4\), we conjecture that
\[
\lim_{N\to\infty}\frac{P_{N,j}}{Q_{N,j}}
=
\frac{\sqrt{10005}}{\pi}.
\]

\section{The Conjectures}

This section contains open problems, presented here for the first time.
In each case, the statement is written in the form of an explicit formula so that it can be checked numerically. 

\subsection{An integral over knot polynomial roots expressing $\pi^2$}
Let $A_{7_{2}}(M, L)$ denote the $A$-polynomial for the prime knot $7_2$. 
\[
\begin{aligned}
A_{7_2}(M,L)
={}& L^5 + L^4\left(
M^{14}-M^{12}+3M^4+4M^2-2
\right) \\
&+ L^3\left(
-2M^{18}+5M^{16}+M^{14}-4M^{12}
+6M^8+5M^6+2M^4-4M^2+1
\right) \\
&+ L^2\left(
M^{22}-4M^{20}+2M^{18}+5M^{16}
+6M^{14}-4M^{10}+M^8+5M^6-2M^4
\right) \\
&+ L\left(
-2M^{22}+4M^{20}+3M^{18}-M^{10}+M^8
\right)+ M^{22}.
\end{aligned}
\]
Let $\alpha \approx 0.349269\ldots$ denote the real solution of $A_{7_{2}}\left( \alpha, \alpha^{1/2} \right) = 0$ closest to $0.349269$.
 Let $\beta \approx 0.406813\ldots$ denote the real solution of $A_{7_{2}}(\beta, \beta) = 0$ closest to $0.406813$
 Let $y = y(x)$ denote the curve satisfying $A(x, y(x)) = 0$ for $\alpha \leq x \leq \beta$ and such that $y(x)$ is positive and decreasing with $y'(x) \leq -2$. 
 Prove the following formula, which was discovered empirically: 
 $$ \frac{4 \pi ^2}{85} = \int_{\alpha}^{\beta} \left( \log x \frac{dy}{y} - \log y \frac{dx}{x} \right). $$

\textbf{Context.}
This result is a surprisingly simple conjectural closed form of a Godbillon–Vey type Knot invariant [\cite{Khoi2008}].

\subsection{Optimality of Ap\'ery's irrationality-measure bound for $\zeta(3)$}

The sequences $(a_n)_{n\ge 0}$ and $(b_n)_{n\ge 0}$ are  solutions of Ap\'ery's recurrence
\[
  (n+1)^3 u_{n+1}
  -\bigl(34n^3+51n^2+27n+5\bigr)u_n
  +n^3 u_{n-1}=0
  \qquad (n\ge 1),
\]
with initial values
\[
  a_0=0,\quad a_1=6,
  \qquad
  b_0=1,\quad b_1=5.
\]
Define $d_n\coloneqq \operatorname{lcm}(1,\dots,n)^3$. Note that $d_na_n,d_nb_n\in \mathbb{Z}$ as shown in [\cite{Apery1979}].

Prove:
\[
\operatorname{gcd}(d_na_n,d_nb_n) = e^{o(n)}.
\]
\textbf{Context.}
The $a_n,b_n$ sequences were generated by Ap\'ery 
in his proof of the irrationality of $\zeta(3)$ [\cite{Apery1979}].
Specifically, considering the linear forms $a_n-\zeta(3) b_n$ 
gives a bound on the irrationality measure of $\zeta(3)$, by
demonstrating that $d_n a_n,d_nb_n\in\mathbb{Z}$. 

The above conjecture shows that the irrationality-measure bound
found by Ap\'ery is tight for his particular formula.

\section{Discussion - proof in the age of AI}\label{section-discussion}

The process of collecting these examples has resurfaced a broader question: what should count as a sufficient proof in the era of AI?
Interactive theorem provers (such as Lean [\cite{moura2015lean}]) are expected to be the gold standard for trust in AI output. In practice, however, many recent AI achievements in mathematics are considered complete without this formalization (e.g., [\cite{openai2026unitdistanceBlog, alon2026remarksdisproofunitdistance}]).
The problems in this challenge sharpen this tension: many are computational in nature, 
and therefore provide a concrete setting in which algorithmic derivations carried out by computer algebra systems can be evaluated as candidates for mathematical proof.

For the purpose of this challenge, 
we consider a derivation carried out using symbolic libraries within established computer algebra systems as sufficient evidence of a valid solution.
This convention should not be interpreted as replacing the strict standards required for a rigorous published proof.
It is motivated by common practice in experimental mathematics, and by the computational nature of the results.
For this reason, problems of the kind presented here may serve as a useful test case for how AI can help create the missing link between numerical discovery, symbolic computation, and rigorous proof.

\printbibliography

@article{Raayoni2021,
   author = {Gal Raayoni and Shahar Gottlieb and Yahel Manor and George Pisha and Yoav Harris and Uri Mendlovic and Doron Haviv and Yaron Hadad and Ido Kaminer},
   doi = {10.1038/s41586-021-03229-4},
   issn = {14764687},
   issue = {7844},
   journal = {Nature},
   title = {Generating conjectures on fundamental constants with the Ramanujan Machine},
   volume = {590},
   year = {2021},
}

@article{Lagarias2013,
   author = {Jeffrey C. Lagarias},
   doi = {10.1090/S0273-0979-2013-01423-X},
   issn = {02730979},
   issue = {4},
   journal = {Bulletin of the American Mathematical Society},
   title = {Euler's constant: Euler's work and modern developments},
   volume = {50},
   year = {2013},
}

@article{elimelech2023algorithm,
  title={Algorithm-assisted discovery of an intrinsic order among mathematical constants},
  author={Elimelech, Rotem and David, Ofir and Mengual, Carlos De la Cruz and Kalisch, Rotem and Berndt, Wolfgang and Shalyt, Michael and Silberstein, Mark and Hadad, Yaron and Kaminer, Ido},
  journal={arXiv preprint arXiv:2308.11829},
  year={2023}
}

@inproceedings{GSM-Symbolic-Numeric-Variation,
title={{GSM}-Symbolic: Understanding the Limitations of Mathematical Reasoning in Large Language Models},
author={Seyed Iman Mirzadeh and Keivan Alizadeh and Hooman Shahrokhi and Oncel Tuzel and Samy Bengio and Mehrdad Farajtabar},
booktitle={The Thirteenth International Conference on Learning Representations},
year={2025},
url={https://openreview.net/forum?id=AjXkRZIvjB}
}

@misc{Humanity-Last_Exam,
      title={Humanity's Last Exam}, 
      author={Long Phan and ...(others)... and Dan Hendrycks},
      year={2025},
      eprint={2501.14249},
      archivePrefix={arXiv},
      primaryClass={cs.LG},
      url={https://arxiv.org/abs/2501.14249}, 
}

@inproceedings{moura2015lean,
  title={The Lean theorem prover (system description)},
  author={de Moura, Leonardo and Kong, Soonho and Avigad, Jeremy and Van Doorn, Floris and von Raumer, Jakob},
  booktitle={International Conference on Automated Deduction},
  pages={378--388},
  year={2015},
  organization={Springer}
}

@article{AbouzaidEtAl2026,
  author = {Mohammed Abouzaid and Andrew J. Blumberg and Martin Hairer and others},
  title = {First Proof},
  journal = {arXiv preprint arXiv:2602.05192},
  year = {2026},
  doi = {10.48550/arXiv.2602.05192},
  url = {https://arxiv.org/abs/2602.05192}
}

@article{GlazerEtAl2024,
  author = {Elliot Glazer and Ege Erdil and Tamay Besiroglu and others},
  title = {FrontierMath: A Benchmark for Evaluating Advanced Mathematical Reasoning in AI},
  journal = {arXiv preprint arXiv:2411.04872},
  year = {2024},
  doi = {10.48550/arXiv.2411.04872},
  url = {https://arxiv.org/abs/2411.04872}
}

@article{ZhangEtAl2025,
  author = {Jie Zhang and Cezara Petrui and Kristina Nikolic and Florian Tramer},
  title = {RealMath: A Continuous Benchmark for Evaluating Language Models on Research-Level Mathematics},
  journal = {arXiv preprint arXiv:2505.12575},
  year = {2025},
  doi = {10.48550/arXiv.2505.12575},
  url = {https://arxiv.org/abs/2505.12575}
}

@article{Ramanujan1914,
  author = {Srinivasa Ramanujan},
  title = {Modular Equations and Approximations to $\pi$},
  journal = {The Quarterly Journal of Pure and Applied Mathematics},
  volume = {45},
  pages = {350--372},
  year = {1914}
}

@incollection{Apery1979,
  author = {Roger Ap\'ery},
  title = {Irrationalit\'e de $\zeta(2)$ et $\zeta(3)$},
  booktitle = {Journ\'ees Arithm\'etiques de Luminy},
  series = {Ast\'erisque},
  volume = {61},
  pages = {11--13},
  publisher = {Soci\'et\'e Math\'ematique de France},
  year = {1979}
}

@article{Beukers1979,
  author = {Frits Beukers},
  title = {A Note on the Irrationality of $\zeta(2)$ and $\zeta(3)$},
  journal = {Bulletin of the London Mathematical Society},
  volume = {11},
  number = {3},
  pages = {268--272},
  year = {1979},
  doi = {10.1112/blms/11.3.268}
}

@article{BallRivoal2001,
  author = {Keith Ball and Tanguy Rivoal},
  title = {Irrationalit\'e d'une infinit\'e de valeurs de la fonction z\^eta aux entiers impairs},
  journal = {Inventiones Mathematicae},
  volume = {146},
  number = {1},
  pages = {193--207},
  year = {2001},
  doi = {10.1007/s002220100168}
}

@article{ChamberlandStraub2021,
  author = {Marc Chamberland and Armin Straub},
  title = {Ap\'ery Limits: Experiments and Proofs},
  journal = {American Mathematical Monthly},
  volume = {128},
  number = {9},
  pages = {811--824},
  year = {2021},
  doi = {10.1080/00029890.2021.1958005}
}

@article{Zudilin2001,
  author = {Wadim Zudilin},
  title = {One of the Numbers $\zeta(5)$, $\zeta(7)$, $\zeta(9)$, $\zeta(11)$ is Irrational},
  journal = {Russian Mathematical Surveys},
  volume = {56},
  number = {4},
  pages = {774--776},
  year = {2001},
  doi = {10.1070/RM2001v056n04ABEH000427}
}

@article{RivoalZudilin2003,
  author = {Tanguy Rivoal and Wadim Zudilin},
  title = {Diophantine Properties of Numbers Related to Catalan's Constant},
  journal = {Mathematische Annalen},
  volume = {326},
  number = {4},
  pages = {705--721},
  year = {2003},
  doi = {10.1007/s00208-003-0420-2}
}

@incollection{KontsevichZagier2001,
  author = {Maxim Kontsevich and Don Zagier},
  title = {Periods},
  booktitle = {Mathematics Unlimited---2001 and Beyond},
  pages = {771--808},
  publisher = {Springer},
  year = {2001}
}

@article{Khoi2008,
  author  = {V. T. Khoi},
  title   = {On the Integral of {$\log x\,\frac{dy}{y}-\log y\,\frac{dx}{x}$} over the {$A$}-Polynomial Curves},
  journal = {Acta Mathematica Vietnamica},
  year    = {2008},
  volume  = {33},
  number  = {3},
  pages   = {519--528},
  url     = {https://arxiv.org/abs/0811.2725},
}

@article{weinbaum2025conservative,
  title={On Conservative Matrix Fields: Continuous Asymptotics and Arithmetic},
  author={Weinbaum, Shachar and Leibtag, Elyasheev and Kalisch, Rotem and Shalyt, Michael and Kaminer, Ido},
  journal={arXiv preprint arXiv:2507.08138},
  year={2025}
}

@misc{ASyMOB,
      title={ASyMOB: Algebraic Symbolic Mathematical Operations Benchmark}, 
      author={Michael Shalyt and Rotem Elimelech and Ido Kaminer},
      year={2025},
      eprint={2505.23851},
      archivePrefix={arXiv},
      primaryClass={cs.CL},
      url={https://arxiv.org/abs/2505.23851}, 
}

@misc{openai2026unitdistanceBlog,
  author       = {{OpenAI}},
  title        = {An OpenAI Model Has Disproved a Central Conjecture in Discrete Geometry},
  year         = {2026},
  month        = may,
  howpublished = {\url{https://openai.com/index/model-disproves-discrete-geometry-conjecture/}},
  note         = {OpenAI Blog, accessed 2026-06-18}
}

@misc{alon2026remarksdisproofunitdistance,
      title={Remarks on the disproof of the unit distance conjecture}, 
      author={Noga Alon and Thomas F. Bloom and W. T. Gowers and Daniel Litt and Will Sawin and Arul Shankar and Jacob Tsimerman and Victor Wang and Melanie Matchett Wood},
      year={2026},
      eprint={2605.20695},
      archivePrefix={arXiv},
      primaryClass={math.CO},
      url={https://arxiv.org/abs/2605.20695}, 
}

@misc{garre2026riemannbenchbenchmarkmoonshotmathematics,
      title={Riemann-Bench: A Benchmark for Moonshot Mathematics}, 
      author={Suhaas Garre and Erik Knutsen and Sushant Mehta and Edwin Chen},
      year={2026},
      eprint={2604.06802},
      archivePrefix={arXiv},
      primaryClass={cs.AI},
      url={https://arxiv.org/abs/2604.06802}, 
}

@misc{peyronnet2026lemmabenchliveresearchlevelbenchmark,
      title={LemmaBench: A Live, Research-Level Benchmark to Evaluate LLM Capabilities in Mathematics}, 
      author={Antoine Peyronnet and Fabian Gloeckle and Amaury Hayat},
      year={2026},
      eprint={2602.24173},
      archivePrefix={arXiv},
      primaryClass={cs.AI},
      url={https://arxiv.org/abs/2602.24173}, 
}

\end{document}